# Quadruple Shehu Transform and Its Applications

D.D. Pawar[1], G.G. Bhuttampalle[1], S.B. Chavhan[2], Wagdi F.S. Ahmed[1], R.D. Kadam[1]
School of Mathematical Sciences,
[1]Swami Ramanand Teerth Marathwada University, Nanded-431606, India
[2]D.B. College, Bhokar, Nanded-431606, India.
E-mail: dypawar@yahoo.com, bgovardhan14@gmail.com and sbcmath2015@gmail.com

**ABSTRACT:** In this current article, we introduce the quadruple Shehu transform and it's inverse. We also introduce some properties of quadruple Shehu transform. The Convolution theorem and its proof are also discussed. Further, to solve homogeneous and non-homogeneous partial differential equation we use this transform.

**INTRODUCTION:**
One of the simplest and most efficient ways to address issues in mathematical physics, applied mathematics, and engineering science that are specified by differential equations, difference equations, and integral equations is to use integral transforms. The basic concept behind the methodology of application is to change the unknown function of a complex variable i.e. t, into another function. This allows the corresponding differential equation to be directly reduced into algebraic equation in the new variable or a differential equation of a lower dimension. Laplace transform is one of the main integral transforms (Papoulis( 1957), Rehman and others(2014), Debnath and Bhatta(2015), Eltayeb and Kilicman (2010), Debnath (2016), Ahmed WFS et al.(2021), Elzaki T. (2012) and Dhunde et al. (2013)), Sumudu transform (Kilicman and Gadain (2010), Mahdy et al. (2015), Mahdy et al. (2015)a, Mechee and Naeemah (2020) and Ahmed WFS and Pawar D.D.(2020)). A transformation of the coordinates and the functions can be used to convert between the Aboodh transform (Aboodh, 2013), Elzaki transform (Elzaki, 2011), Variational homotopy perturbation method (Mahdy et al., 2015b), Alternative variational iteration method (Mtawal et al., 2020), and one from the other. For the purpose of solving differential equations in the time domain, Maitama and Zhao (2019) recently developed a new type of integral transform in 2019 as an extension of both the Laplace transform and the Sumudu transform. They also provided several theorems on this transform.

The results of Abdon (2013), Bokhari A et al. (2019), Thakur et al. (2018), Alfaqeih (2020) and Sameehah et al. (2021) are extended and generalized in this study; specifically, the one-dimensional Shehu transform is extended into the four-dimensional Shehu transform (quadruple Shehu transform) and several examples are given to demonstrate the usefulness of our findings.

**PRELIMINARIES**
In current section we recall definitions of Shehu, Double Shehu and Triple Shehu transform.
**Definition 1** The Shehu transform (Maitama and Zhao, 2019) is defined over the set of functions

$$A=\left\{f(q):\exists M,\mu_1,\mu_2>0,|f(q)|<Me^{\frac{|q|}{\mu_j}},\text{if } q\in(-1)^j*[0,\infty)\right\}$$

By the given formula

$$\mathbb{H}(f(q))=F(h,m)=\int_0^\infty e^{-\left(\frac{hq}{m}\right)}f(q)dq \qquad (1)$$

$$=\lim_{\alpha\to\infty}\int_0^\infty e^{-\left(\frac{hq}{m}\right)}f(q)dq;\ \ h,m>0$$

And its inverse is defined by the following formula

$$Ƕ^{-1}(f(q)) = [F(h,m)] = f(q), q \geq 0,$$

$$= \frac{1}{2\pi i} \int_{\alpha - i\infty}^{\alpha + i\infty} \frac{1}{m} e^{\left(\frac{hq}{m}\right)} F(h,m) dh \qquad (2)$$

**Definition 2** A real function f(x), $x > 0$, is considered to be in the space $C_m, m \in R$, if their exist a real number $\sigma > m$, so that $f(t) = t^{\sigma} g(x)$, where $g(t) \in [0, \infty)$, and it is said to be in the space $C_{\sigma}^{m}$, if $f^m \in C_{\sigma}, m \in N$. (Podlubny, 1999; He, 2014).

**Definition 3** The left-sided Riemann-Liouville fractional integral of order $\alpha \geq 0$, of a function $f \in C_{\sigma}, \sigma \geq -1$, (Podlubny, 1999; He, 2014) is defined as:

$$J^{\alpha} f(t) = \frac{1}{\Gamma(\alpha)} \int_0^t (t-\tau)^{\alpha-1} f(\tau) d\tau, t, \alpha > 0. \qquad (3)$$

Here $\Gamma(.)$ is a gamma function.

**Definition 4** If $f \in C_m^n, n \in N \cup \{0\}$. The left Caputo fractional derivative if $f$ in the Caputo sense (Podlubny, 1999; He, 2014) is defined as follows:

$$D_t^{\alpha} f(t) = \begin{cases} \frac{1}{\Gamma(n-\alpha)} \int_0^t (t-\tau)^{n-\alpha-1} f^{(n)}(\tau) d\tau & , n-1 < \alpha \leq n, \\ D_t^{\alpha} f(t), & , \alpha = n \end{cases} \qquad (4)$$

**Definition 5** The Mittag-Leffler function $E_{\gamma}(z), E_{\gamma,\beta}(z)$ (Kilbas et al., 2004) are defined as

$$E_{\gamma}(q) = \sum_{r=0}^{\infty} \frac{q^r}{\Gamma(\gamma r + 1)} \qquad (\gamma, q \in C, R(\gamma) > 0) \qquad (5)$$

$$E_{\gamma,\beta}(q) = \sum_{r=0}^{\infty} \frac{q^r}{\Gamma(\gamma r + \beta)} \qquad (\beta, \gamma, q \in C, R(\gamma) > 0, R(\beta) > 0) \qquad (6)$$

These functions are generalization of the exponential function. Some special cases of the Mittag-Leffler function are as follows:

$$E_1(t) = e^t, E_{(\gamma,1)}(t) = E_{\gamma}(t)$$

**Definition 6** The single Shehu transform (ǶT) of a real valued function $f(q,r,s,t)$ with respect to the variables $q, r, s, and\ t$ respectively, (Maitma and Zhao, 2019) are defined by:

$$Ƕ_q(f(q,r,s,t)) = \int_0^{\infty} e^{-\left(\frac{hq}{m}\right)} f(q,r,s,t) dq, \qquad (7)$$

$$Ƕ_r(f(q,r,s,t)) = \int_0^{\infty} e^{-\left(\frac{jr}{n}\right)} f(q,r,s,t) dr, \qquad (8)$$

$$Ƕ_s(f(q,r,s,t)) = \int_0^{\infty} e^{-\left(\frac{ks}{o}\right)} f(q,r,s,t) ds, \qquad (9)$$

$$Ƕ_t(f(q,r,s,t)) = \int_0^{\infty} e^{-\left(\frac{lt}{p}\right)} f(q,r,s,t) dt, \qquad (10)$$

**Definition 7** The double Shehu transform $Ƕ^2$ as a function f(q, r) (Alfaqeih and Misirli, 2020) is defined over the Set of functions

$$A = \left\{ f(q,r): \exists M, \mu_1, \mu_2 > 0, |f(q,r)| < Me^{\frac{|q+r|}{\mu_j}}, \text{if } (q,r) \in (R)^2, j = 1,2 \right\}$$

By the following formula

The double Shehu transform of the function $f(q,r)$ is given by,

$$Ⴆ_{qr}^{2}\big(f(q,r)\big) = F[(h,j),(m,n)] = \int_0^\infty \int_0^\infty e^{-\left(\frac{hq}{m}+\frac{jr}{n}\right)} f(q,r)\,dq\,dr, \tag{11}$$

And the inverse double Shehu transform is defined by

$$Ⴆ_{qr}^{-2}F[(h,j),(m,n)] = f(q,r)$$
$$= \frac{1}{2\pi i}\int_{\alpha-i\infty}^{\alpha+i\infty} \frac{1}{h} e^{\left(\frac{hq}{m}\right)} \left[ \frac{1}{2\pi i}\int_{\beta-i\infty}^{\beta+i\infty} \frac{1}{k} e^{\left(\frac{jr}{n}\right)} H_2\big(f(q,r)\big)\,dj \right] dh \tag{12}$$

**Definition 8** The double Shehu transform $Ⴆ^2$ of a function $f(q,r,s,t)$with respect to $qr, qs, qt, rs, rt \text{ and } st$ respecctively (Alfaqeih and Misirli, 2020) are defined by:

$$Ⴆ_{qr}^{2}\big(f(q,r,s,t)\big) = \int_0^\infty \int_0^\infty e^{-\left(\frac{hq}{m}+\frac{jr}{n}\right)} f(q,r,s,t)\,dq\,dr, \tag{13}$$

$$Ⴆ_{rs}^{2}\big(f(q,r,s,t)\big) = \int_0^\infty \int_0^\infty e^{-\left(\frac{jr}{n}+\frac{ks}{o}\right)} f(q,r,s,t)\,dr\,ds,$$

$$Ⴆ_{st}^{2}\big(f(q,r,s,t)\big) = \int_0^\infty \int_0^\infty e^{-\left(\frac{ks}{o}+\frac{lt}{p}\right)} f(q,r,s,t)\,ds\,dt,$$

$$Ⴆ_{qs}^{2}\big(f(q,r,s,t)\big) = \int_0^\infty \int_0^\infty e^{-\left(\frac{hq}{m}+\frac{ks}{o}\right)} f(q,r,s,t)\,dq\,ds,$$

$$Ⴆ_{qt}^{2}\big(f(q,r,s,t)\big) = \int_0^\infty \int_0^\infty e^{-\left(\frac{hq}{m}+\frac{lt}{p}\right)} f(q,r,s,t)\,dq\,dt,$$

$$Ⴆ_{rt}^{2}\big(f(q,r,s,t)\big) = \int_0^\infty \int_0^\infty e^{-\left(\frac{jr}{n}+\frac{lt}{p}\right)} f(q,r,s,t)\,dr\,dt,$$

**Definition 9** The triple Shehu transform $Ⴆ^3$ of a function $f(q,r,s)$ (Sameehah R.A. et al,2021) is defined over the set of functions

$$A = \left\{ f(q,r,s): \exists M, \mu_1, \mu_2 > 0, |f(q,r,s)| < Me^{\frac{|q+r+s|}{\mu_j}}, \text{if } (q,r,s) \in (R)^2, j = 1,2,3 \right\}$$

By the following formula

The triple Shehu transform of the function $f(q,r,s)$ is given by

$$\text{Ƕ}^3_{qrs}\big(f(q,r,s)\big) = F[(h,j,k),(m,n,o)]$$

$$= \int_0^\infty \int_0^\infty \int_0^\infty e^{-\left(\frac{hq}{m}+\frac{jr}{n}+\frac{ks}{o}\right)} f(q,r,s)\,dq\,dr\,ds, \qquad (14)$$

And the inverse triple Shehu transform is defined by

$$\text{Ƕ}^3_{qrs} F[(i,j,k),(m,n,o)] = f(q,r,s)$$

$$= \frac{1}{2\pi i}\int_{\alpha-i\infty}^{\alpha+i\infty} \frac{1}{m} e^{\left(\frac{hq}{m}\right)} \left[\frac{1}{2\pi i}\int_{\beta-i\infty}^{\beta+i\infty} \frac{1}{k} e^{\left(\frac{jr}{n}\right)} \left(\frac{1}{2\pi i}\int_{\gamma-i\infty}^{\gamma+i\infty} \frac{1}{l} e^{\left(\frac{ks}{o}\right)} H_3(f(q,r,s))\,dk\right) dj\right] dh \qquad (15)$$

**Definition 10** The triple Shehu transform $\text{Ƕ}^3$ of a function $f(q,r,s,t)$ (Sameehah R.A. et al, 2021) with respect to $qrs, qrt, qst, rst, rtq, rsq, tqr, trs$ *and* $tqs$ respectively, are defined by:

$$\text{Ƕ}^3_{qrs}\big(f(q,r,s,t)\big) = \int_0^\infty \int_0^\infty \int_0^\infty e^{-\left(\frac{hq}{m}+\frac{jr}{n}+\frac{ks}{o}\right)} f(q,r,s,t)\,dq\,dr\,ds, \qquad (16)$$

$$\text{Ƕ}^3_{qrt}\big(f(q,r,s,t)\big) = \int_0^\infty \int_0^\infty \int_0^\infty e^{-\left(\frac{hq}{m}+\frac{jr}{n}+\frac{lt}{p}\right)} f(q,r,s,t)\,dq\,dr\,dt,$$

$$\text{Ƕ}^3_{qst}\big(f(q,r,s,t)\big) = \int_0^\infty \int_0^\infty \int_0^\infty e^{-\left(\frac{hq}{m}+\frac{ks}{o}+\frac{lt}{p}\right)} f(q,r,s,t)\,dq\,ds\,dt,$$

$$\text{Ƕ}^3_{rst}\big(f(q,r,s,t)\big) = \int_0^\infty \int_0^\infty \int_0^\infty e^{-\left(\frac{jr}{n}+\frac{ks}{o}+\frac{lt}{p}\right)} f(q,r,s,t)\,dr\,ds\,dt,$$

$$\text{Ƕ}^3_{rtq}\big(f(q,r,s,t)\big) = \int_0^\infty \int_0^\infty \int_0^\infty e^{-\left(\frac{jr}{n}+\frac{lt}{p}+\frac{hq}{m}\right)} f(q,r,s,t)\,dr\,dt\,dq,$$

$$\text{Ƕ}^3_{rsq}\big(f(q,r,s,t)\big) = \int_0^\infty \int_0^\infty \int_0^\infty e^{-\left(\frac{jr}{n}+\frac{ks}{o}+\frac{hq}{m}\right)} f(q,r,s,t)\,dr\,ds\,dq,$$

$$\text{Ƕ}^3_{tqr}\big(f(q,r,s,t)\big) = \int_0^\infty \int_0^\infty \int_0^\infty e^{-\left(\frac{lt}{p}+\frac{hq}{m}+\frac{jr}{n}\right)} f(q,r,s,t)\,dt\,dq\,dr,$$

$$\text{Ƕ}^3_{trs}\big(f(q,r,s,t)\big) = \int_0^\infty \int_0^\infty \int_0^\infty e^{-\left(\frac{lt}{p}+\frac{jr}{n}+\frac{ks}{o}\right)} f(q,r,s,t)\,dt\,ds\,dr,$$

$$\text{Ƕ}^3_{tqs}\big(f(q,r,s,t)\big) = \int_0^\infty \int_0^\infty \int_0^\infty e^{-\left(\frac{lt}{p}+\frac{hq}{m}+\frac{ks}{o}\right)} f(q,r,s,t)\,dt\,dq\,ds;$$

## RESULTS

Here we are introduce the definition of fourth Shehu transform and fourth Shehu transform of partial and fractional derivatives which are used further in this paper, moreover we apply fourth Shehu transform for some basic functions.

**Definition 11** Let f is a continuous function of four variables. Then the fourth Shehu transform $Ѣ^4$ of $f(q,r,s,t)$ is defined by

$$Ѣ^4_{qrst}\big(f(q,r,s,t)\big) = F[(h,j,k,l),(m,n,o,p)]$$

$$= \int_0^\infty \int_0^\infty \int_0^\infty \int_0^\infty e^{-\left(\frac{hq}{m}+\frac{jr}{n}+\frac{ks}{o}+\frac{lt}{p}\right)} f(q,r,s,t)dqdrdsdt, \qquad (17)$$

Where $q,r,s,t \geq 0$ *and* $h,j,k,l,m,n,o,p$ are Shehu variables, provided the integral exist and the inverse fourth Shehu transform is defined by

$$Ѣ^4_{qrst}F[(h,j,k,l),(m,n,o,p)] = f(q,r,s,t)$$

$$= \frac{1}{2\pi i}\int_{\alpha-i\infty}^{\alpha+i\infty} \frac{1}{m} e^{\left(\frac{hq}{m}\right)} \left[\frac{1}{2\pi i}\int_{\beta-i\infty}^{\beta+i\infty} \frac{1}{n} e^{\left(\frac{jr}{n}\right)} \left(\frac{1}{2\pi i}\int_{\gamma-i\infty}^{\gamma+i\infty} \frac{1}{o} e^{\left(\frac{ks}{o}\right)} \left(\frac{1}{2\pi i}\int_{\delta-i\infty}^{\delta+i\infty} \frac{1}{p} e^{\left(\frac{lt}{p}\right)} H_4(f(x,y,z,t))dl\right)dk\right)dj\right]dh \quad (18)$$

## Existence and uniqueness of quadruple Shehu transform

Now, we define existence and uniqueness quadruple Shehu transform and prove it.

## Definition 12

A function $f(q,r,s,t)$ is said to be of exponential order $h>0, j>0, k>0, l>0$ as $q,r,s,t \to \infty$ if there are positive constants $M, X, Y, Z$ *and* $T$ (Alfaqeih, 2019) such that
$|f(q,r,s,t)| \leq Me^{(hq+jr+ks+lt)}$ for all $q>X, r>Y, s>Z, t>T$ and we write,
$f(q,r,s,t) = O(e^{(hq+jr+ks+lt)})$ $(as\ q,r,s,t \to \infty)$
Or equivalently,

$$\sup_{q,r,s,t>0}\left[\frac{|f(q,r,s,t)|}{e^{(aq+br+cs+dt)}}\right] < \infty$$

**Theorem 1**

Let $f(q,r,s,t)$ be a continuous function on the interval $(0,X),(0,Y),(0,Z),(0,T)$ and of exponential order $e^{(aq+br+cs+dt)}$ then the Fourth Shehu transform of $f(q,r,s,t)$ exists.

$$\forall h > am, j > bn, k > co, l < dp$$

$\overset{Proof}{\Longrightarrow}$ Let $f(q,r,s,t)$ be of exponential order $e^{(aq+br+cs+dt)}$ such that
$|f(q,r,s,t)| \leq Me^{(aq+br+cs+dt)}$ for all $q>X, r>Y, s>Z, t>T$
Then we have

$$\left|Ѣ^4_{qrst}(f(q,r,s,t))\right| = \left|\int_0^\infty \int_0^\infty \int_0^\infty \int_0^\infty e^{-\left(\frac{hq}{m}+\frac{jr}{n}+\frac{ks}{o}+\frac{lt}{p}\right)} f(q,r,s,t)dqdrdsdt\right|$$

$$\leq \int_0^\infty \int_0^\infty \int_0^\infty \int_0^\infty e^{-\left(\frac{hq}{m}+\frac{jr}{n}+\frac{ks}{o}+\frac{lt}{p}\right)} |f(q,r,s,t)|dqdrdsdt$$

$$\leq \int_0^\infty \int_0^\infty \int_0^\infty \int_0^\infty e^{-\left(\frac{hq}{m}+\frac{jr}{n}+\frac{ks}{o}+\frac{lt}{p}\right)} e^{(aq+br+cs+dt)}dqdrdsdt$$

$$= M\int_0^\infty e^{-\left(\frac{l-dp}{p}\right)t}\int_0^\infty e^{-\left(\frac{k-co}{o}\right)z}\int_0^\infty e^{-\left(\frac{j-bn}{n}\right)y}\int_0^\infty e^{-\left(\frac{h-am}{m}\right)x} dqdrdsdt$$

$$= \frac{M\, mnop}{(h-am)(j-bn)(k-co)(l-dp)}$$

Thus the proof is complete.

In the next theorem, we show that $f(q,r,s,t)$ can be uniquely obtained.

**Theorem 2**

Let $F_1[(h,j,k,l),(m,n,o,p)]$ and $F_2[(h,j,k,l),(m,n,o,p)]$ be the Shehu transform of the continuous functions $f_1(q,r,s,t)$ & $f_2(q,r,s,t)$ defined for $q,r,s,t \geq 0$ resp. If $F_1[(h,j,k,l),(m,n,o,p)] = F_2[(h,j,k,l),(m,n,o,p)]$ then $f_1(q,r,s,t) = f_2(q,r,s,t)$

$\overset{Proof}{\Longrightarrow}$ If we presume $\alpha,\beta,\gamma,\delta$ to be sufficiently large, then since

$$f(q,r,s,t) = \frac{1}{2\pi i}\int_{\alpha-i\infty}^{\alpha+i\infty}\frac{1}{m}e^{\left(\frac{hq}{m}\right)}\left\{\frac{1}{2\pi i}\int_{\beta-i\infty}^{\beta+i\infty}\frac{1}{n}e^{\left(\frac{jr}{n}\right)}\left[\frac{1}{2\pi i}\int_{\gamma-i\infty}^{\gamma+i\infty}\frac{1}{o}e^{\left(\frac{ks}{o}\right)}\left(\frac{1}{2\pi i}\int_{\delta-i\infty}^{\delta+i\infty}\frac{1}{p}e^{\left(\frac{lt}{p}\right)}F[(h,j,k,l),(m,n,o,p)]dl\right)dk\right]dj\right\}dh$$

We deduce that

$$f_1(q,r,s,t) = \frac{1}{2\pi i}\int_{\alpha-i\infty}^{\alpha+i\infty}\frac{1}{m}e^{\left(\frac{hq}{m}\right)}\left\{\frac{1}{2\pi i}\int_{\beta-i\infty}^{\beta+i\infty}\frac{1}{n}e^{\left(\frac{jr}{n}\right)}\left[\frac{1}{2\pi i}\int_{\gamma-i\infty}^{\gamma+i\infty}\frac{1}{o}e^{\left(\frac{ks}{o}\right)}\left(\frac{1}{2\pi i}\int_{\delta-i\infty}^{\delta+i\infty}\frac{1}{p}e^{\left(\frac{lt}{p}\right)}F_1[(s,q,r,p),(u,v,k,l)]dl\right)dk\right]dj\right\}dh$$

$$f(q,r,s,t) = \frac{1}{2\pi i}\int_{\alpha-i\infty}^{\alpha+i\infty}\frac{1}{m}e^{\left(\frac{hq}{m}\right)}\left\{\frac{1}{2\pi i}\int_{\beta-i\infty}^{\beta+i\infty}\frac{1}{n}e^{\left(\frac{jr}{n}\right)}\left[\frac{1}{2\pi i}\int_{\gamma-i\infty}^{\gamma+i\infty}\frac{1}{o}e^{\left(\frac{ks}{o}\right)}\left(\frac{1}{2\pi i}\int_{\delta-i\infty}^{\delta+i\infty}\frac{1}{p}e^{\left(\frac{lt}{p}\right)}F_2[(s,q,r,p),(u,v,k,l)]dl\right)dk\right]dj\right\}dh$$

$$= f_2(q,r,s,t)$$

Thus the uniqueness of quadruple Shehu transform is proved.

## Fourth Shehu transform of some elementary functions

(1) If $f(q,r,s,t) = A; q,r,s,t > 0$ then $Ƕ^4_{qrst}(A) = A\frac{mnop}{hjkl}$

$$Ƕ^4_{qrst}[f(q,r,s,t)] = \int_0^\infty\int_0^\infty\int_0^\infty\int_0^\infty e^{-\left(\frac{hq}{m}+\frac{jr}{n}+\frac{ks}{o}+\frac{lt}{p}\right)}f(q,r,s,t)dxdydzdt$$

(2) If $f(q,r,s,t)qrst$. Then $Ƕ^4_{qrst}(q\ r\ s\ t) = \left(\frac{mnop}{hjkl}\right)^2$

(3) If $f(q,r,s,t) = e^{aq+br+cs+dt}$ then $Ƕ^4_{qrst}\left(e^{aq+br+cs+dt}\right) = \frac{mnop}{(h-am)(j-bn)(k-co)(l-dp)}$

(4) If $f(q,r,s,t) = e^{i(aq+br+cs+dt)}$ then

$$Ƕ^4_{qrst}\left(e^{i(aq+br+cs+dt)}\right) = \frac{mnop}{(h-iam)(j-ibn)(k-ico)(l-idp)}$$

$$= \frac{mnop\begin{Bmatrix}hjkl - hjocdp - ambnkl + ambnocdp - hbnkdp - hbnlop - ajmkdp - amjloc + \\ i[hjkdp + hjloc - ambnkdp - ambnloc + klhbn + klamj - hbnocdp - amjpcdp]\end{Bmatrix}}{(h^2+a^2m^2)(j^2+b^2n^2)(k^2+o^2c^2)(l^2+d^2p^2)}$$

Consequently,

$$Ƕ^4_{qrst}[cos(aq+br+cs+dt)] = \frac{mnop[hjkl - hjocdp - ambnkl + ambnocdp - hbnkdp - hbnlop - ajmkdp - amjloc]}{(h^2+a^2m^2)(j^2+b^2n^2)(k^2+o^2c^2)(l^2+d^2p^2)}$$

$$Ƕ^4_{qrst}[sin(aq+br+cs+dt)] = \frac{hjkdp+hjloc-ambnkdp-ambnloc+klhbn+klamj-hbnocdp-amjpcdp}{(h^2+a^2m^2)(j^2+b^2n^2)(k^2+o^2c^2)(l^2+d^2p^2)} \quad (19)$$

If $f(q,r,s,t) = f_1(q)f_2(r)f_3(s)f_4(t)$ then

$$Ƕ^4_{qrst}[f(q,r,s,t)] = \int_0^\infty\int_0^\infty\int_0^\infty\int_0^\infty e^{-\left(\frac{hq}{m}+\frac{jr}{n}+\frac{ks}{o}+\frac{lt}{p}\right)} f_1(q)f_2(r)f_3(s)f_4(t)dqdrdsdt$$

$$= \int_0^\infty e^{\left(\frac{-hq}{m}\right)} f_1(q)dq \int_0^\infty e^{\left(\frac{-jr}{n}\right)} f_2(r)dr \int_0^\infty e^{\left(\frac{-ks}{o}\right)} f_3(s)ds \int_0^\infty e^{\left(\frac{-lt}{p}\right)} f_4(t)dt$$

$$= Ƕ_q f_1(q) Ƕ_r f_2(r) Ƕ_s f_3(s) Ƕ_t f_4(t)$$

So that

$$Ƕ^4(cosq\ cosr\ coss\ cost) = \left(\frac{hm}{(h^2+m^2)}\right)\left(\frac{jn}{(j^2+n^2)}\right)\left(\frac{ko}{(k^2+o^2)}\right)\left(\frac{lp}{(l^2+p^2)}\right)$$

$$Ƕ^4(sinq\ sinr\ sins\ sint) = \left(\frac{m^2}{(h^2+m^2)}\right)\left(\frac{n^2}{(j^2+n^2)}\right)\left(\frac{o^2}{(k^2+o^2)}\right)\left(\frac{p^2}{(l^2+p^2)}\right)$$

$$Ƕ^4_{qrst}(q^{n_1}r^{n_2}s^{n_3}t^{n_4}) = n_1!\left(\frac{m}{h}\right)^{n_1+1} n_2!\left(\frac{n}{j}\right)^{n_2+1} n_3!\left(\frac{o}{k}\right)^{n_3+1} n_4!\left(\frac{p}{l}\right)^{n_4+1}, n_1,n_2,n_3,n_4 = 0,1,2..$$

$$Ƕ^4_{qrst}(qrst)^n = (n!)^4\left(\frac{mnop}{hjkl}\right)^{n+1}, n = 0,1,2,..$$

OR

$$Ƕ^4_{qrst}(q^{n_1}r^{n_2}s^{n_3}t^{n_4}) = \Gamma(n_1+1)\left(\frac{m}{h}\right)^{n_1+1}\Gamma(n_2+1)\left(\frac{n}{j}\right)^{n_2+1}\Gamma(n_3+1)\left(\frac{o}{k}\right)^{n_3+1}\Gamma(n_4+1)\left(\frac{p}{l}\right)^{n_1+1}, n_1,n_2,n_3,n_4 \geq -1$$

**Some main properties of fourth Shehu transform**

**Linearity property:-**

$$Ƕ^4_{qrst}[f(q,r,s,t)] = F[(h,j,k,l),(m,n,o,p)]$$

$$Ƕ^4_{qrst}[g(q,r,s,t)] = F[(h,j,k,l),(m,n,o,p)]$$

Then for any constant $\alpha, \beta$ we have,

$$Ƕ^4_{qrst}[\alpha f(q,r,s,t)+\beta g(q,r,s,t)] = \alpha Ƕ^4_{qrst}[f(q,r,s,t)] + \beta Ƕ^4_{qrst}[f(q,r,s,t)]$$

$\xRightarrow{Proof}$ Using the definition of fourth Shehu transform,

$$Ƕ^4_{qrst}[\alpha f(q,r,s,t)+\beta g(q,r,s,t)]$$

$$= \int_0^\infty\int_0^\infty\int_0^\infty\int_0^\infty e^{-\left(\frac{hq}{m}+\frac{jr}{n}+\frac{ks}{o}+\frac{lt}{p}\right)}[\alpha f(q,r,s,t)+\beta g(q,r,s,t)]dqdrdsdt$$

$$= \alpha\int_0^\infty\int_0^\infty\int_0^\infty\int_0^\infty e^{-\left(\frac{hq}{m}+\frac{jr}{n}+\frac{ks}{o}+\frac{lt}{p}\right)} f(q,r,s,t)dqdrdsdt$$

$$+\beta\int_0^\infty\int_0^\infty\int_0^\infty\int_0^\infty e^{-\left(\frac{hq}{m}+\frac{jr}{n}+\frac{ks}{o}+\frac{lt}{p}\right)} g(q,r,s,t)dqdrdsdt$$

$$= \alpha Ƕ^4_{qrst}[f(q,r,s,t)] + \beta Ƕ^4_{qrst}[f(q,r,s,t)]$$

## Change of scale property:-

Let $f(q,r,s,t)$ be a functions such that

$$Ƕ^4_{qrst}[f(q,r,s,t)] = F[(h,j,k,l),(m,n,o,p)]$$

Then for $a,b,c,d > 0$ we have

$$Ƕ^4_{qrst}[f(aq,br,cs,dt)] = \frac{1}{abcd}F\left[\left(\frac{h}{a},\frac{j}{b},\frac{k}{c},\frac{l}{d}\right),(m,n,o,p)\right]$$

$\xRightarrow{Proof}$ We have,

$$\text{Ꞟ}^4_{qrst}[f(aq,br,cs,dt)] = \int_0^\infty \int_0^\infty \int_0^\infty \int_0^\infty e^{-\left(\frac{hq}{m}+\frac{jr}{n}+\frac{ks}{o}+\frac{lt}{p}\right)} f(aq,br,cs,dt)\,dqdrdsdt$$

Let $w_1 = aq, w_2 = br, w_3 = cs, w_4 = dt$ then

$$\text{Ꞟ}^4_{qrst}[f(w_1,w_2,w_3,w_4)] = \frac{1}{abcd}\int_0^\infty \int_0^\infty \int_0^\infty \int_0^\infty e^{-\left[\frac{\frac{hw_1}{a}}{m}+\frac{\frac{jw_2}{b}}{n}+\frac{\frac{kw_3}{c}}{o}+\frac{\frac{lw_4}{d}}{p}\right]} f(w_1,w_2,w_3,w_4)\,dw_1dw_2dw_3dw_4$$

$$= \frac{1}{abcd}F\left[\left(\frac{h}{a},\frac{j}{b},\frac{k}{c},\frac{l}{d}\right),(m,n,o,p)\right]$$

**First shifting property:-**

Let $f(q,r,s,t)$ be a functions such that

$$\text{Ꞟ}^4_{qrst}[f(q,r,s,t)] = F[(h,j,k,l),(m,n,o,p)]$$

Then for $a,b,c,d > 0$ we have

$$\text{Ꞟ}^4_{qrst}\left[e^{(aq+br+cs+dt)}f(q,r,s,t)\right] = F[(h-am,j-bn,k-co,l-dp),(m,n,o,p)]$$

$\overset{Proof}{\Longrightarrow}$ We have, by definition

$$\text{Ꞟ}^4_{qrst}\left[e^{(aq+br+cs+dt)}f(q,r,s,t)\right] = \int_0^\infty \int_0^\infty \int_0^\infty \int_0^\infty e^{-\left(\frac{hq}{m}+\frac{jr}{n}+\frac{ks}{o}+\frac{lt}{p}\right)} e^{(aq+br+cs+dt)} z f(q,r,s,t)\,dqdrdsdt$$

$$= \int_0^\infty \int_0^\infty \int_0^\infty \int e^{-\left[\left(\frac{(h-am)}{m}\right)x+\left(\frac{(j-bn)}{n}\right)y+\left(\frac{(k-co)}{o}\right)z+\left(\frac{(l-dp)}{p}\right)t\right]} f(q,r,s,t)\,dqdrdsdt$$

$$= F[(h-am,j-bn,k-co,l-dp),(m,n,o,p)]$$

## Fourth Shehu transform of derivative of a function of four variables

1) Fourth Shehu transform of mixed derivative of a function of four variables is given by:

$$\begin{aligned}\text{Ꞟ}^4_{qrst}\left[\frac{\partial^4 f(q,r,s,t)}{\partial q\,\partial r\,\partial s\,\partial t}\right] &= \left(\frac{hjkl}{mnop}\right)F[(h,j,k,l),(m,n,o,p)] - \left(\frac{hjk}{mno}\right)\text{Ꞟ}^3_{qrs}[f(q,r,s,0)] - \left(\frac{hjl}{mnp}\right)\text{Ꞟ}^3_{qrt}[f(q,r,0,t)]\\ &- \left(\frac{hkl}{mop}\right)\text{Ꞟ}^3_{qst}[f(q,0,s,t)] - \left(\frac{jkl}{nop}\right)\text{Ꞟ}^3_{rst}[f(0,r,s,t)] + \left(\frac{hj}{mn}\right)\text{Ꞟ}^2_{qr}[f(q,r,0,0)]\\ &+ \left(\frac{hk}{mo}\right)\text{Ꞟ}^2_{qs}[f(q,0,s,0)] + \left(\frac{hl}{mp}\right)\text{Ꞟ}^2_{qt}[f(q,0,0,t)] + \left(\frac{jk}{no}\right)\text{Ꞟ}^2_{rs}[f(0,r,s,0)] + \left(\frac{jl}{np}\right)\text{Ꞟ}^2_{rt}[f(0,r,0,t)]\\ &+ \left(\frac{kl}{op}\right)\text{Ꞟ}^2_{st}[f(0,0,s,t)] - \left(\frac{h}{m}\right)\text{Ꞟ}_q[f(q,0,0,0)] - \left(\frac{j}{n}\right)\text{Ꞟ}_r[f(0,r,0,0)] - \left(\frac{k}{o}\right)\text{Ꞟ}_s[f(0,0,s,0)]\\ &- \left(\frac{l}{p}\right)\text{Ꞟ}_t[f(0,0,0,t)] + f(0,0,0,0)\end{aligned}$$

2) The Fourth Shehu transform of $n^{th}$ partial derivative of a function of four variables is given by :

$$\text{Ꞟ}^4_{qrst}\left[\frac{\partial^n f(q,r,s,t)}{\partial q^n}\right] = \left(\frac{h}{m}\right)^n F[(h,j,k,l),(m,n,o,p)] - \sum_{m=0}^{n-1}\left(\frac{h}{m}\right)^{n-m-1}\text{Ꞟ}^3_{rst}\left(\frac{\partial^m f(0,r,s,t)}{\partial q^m}\right)$$

$$\text{Ꞟ}^4_{qrst}\left[\frac{\partial^n f(q,r,s,t)}{\partial r^n}\right] = \left(\frac{j}{n}\right)^n F[(h,j,k,l),(m,n,o,p)] - \sum_{m=0}^{n-1}\left(\frac{j}{n}\right)^{n-m-1}\text{Ꞟ}^4_{qst}\left(\frac{\partial^m f(q,0,s,t)}{\partial r^m}\right)$$

$$\text{Ꞟ}^4_{qrst}\left[\frac{\partial^n f(q,r,s,t)}{\partial s^n}\right] = \left(\frac{k}{o}\right)^n F[(h,j,k,l),(m,n,o,p)] - \sum_{m=0}^{n-1}\left(\frac{k}{o}\right)^{n-m-1}\text{Ꞟ}^4_{qrt}\left(\frac{\partial^m f(q,r,0,t)}{\partial s^m}\right)$$

$$\text{Ꞟ}^4_{qrst}\left[\frac{\partial^n f(q,r,s,t)}{\partial t^n}\right] = \left(\frac{l}{p}\right)^n F[(h,j,k,l),(m,n,o,p)] - \sum_{m=0}^{n-1}\left(\frac{l}{p}\right)^{n-m-1}\text{Ꞟ}^3_{qrs}\left(\frac{\partial^m f(q,r,s,0)}{\partial t^m}\right)$$

3) The Fourth Shehu transform of the partial fractional Caputo derivatives of a function of four variables is given by:

$$\text{Ꞟ}^4_{qrst}\left[\frac{\partial^\alpha f(q,r,s,t)}{\partial q^\alpha}\right] = \left(\frac{h}{m}\right)^\alpha F[(h,j,k,l),(m,n,o,p)] - \sum_{m=0}^{n-1}\left(\frac{h}{m}\right)^{\alpha-m-1}\text{Ꞟ}^3_{rst}\left(\frac{\partial^m f(0,r,s,t)}{\partial q^m}\right)$$

$$\text{Ꞟ}^4_{qrst}\left[\frac{\partial^\alpha f(q,r,s,t)}{\partial r^\alpha}\right] = \left(\frac{j}{n}\right)^\alpha F[(h,j,k,l),(m,n,o,p)] - \sum_{m=0}^{n-1}\left(\frac{j}{n}\right)^{\alpha-m-1}\text{Ꞟ}^4_{qst}\left(\frac{\partial^m f(q,0,s,t)}{\partial r^m}\right)$$

$$\text{Ƕ}^4_{\mathrm{qrst}}\left[\frac{\partial^{\alpha}\mathrm{f(q,r,s,t)}}{\partial \mathrm{s}^{\alpha}}\right]=\left(\frac{\mathrm{k}}{\mathrm{o}}\right)^{\alpha}\mathrm{F[(h,j,k,l),(m,n,o,p)]}-\sum_{\mathrm{m}=0}^{\mathrm{n}-1}\left(\frac{\mathrm{k}}{\mathrm{o}}\right)^{\alpha-\mathrm{m}-1}\text{Ƕ}^4_{\mathrm{qrt}}\left(\frac{\partial^{\mathrm{m}}\mathrm{f(q,r,0,t)}}{\partial \mathrm{s}^{\mathrm{m}}}\right)$$

$$\text{Ƕ}^4_{\mathrm{qrst}}\left[\frac{\partial^{\alpha}\mathrm{f(q,r,s,t)}}{\partial \mathrm{t}^{\alpha}}\right]=\left(\frac{\mathrm{l}}{\mathrm{p}}\right)^{\alpha}\mathrm{F[(h,j,k,l),(m,n,o,p)]}-\sum_{\mathrm{m}=0}^{\mathrm{n}-1}\left(\frac{\mathrm{l}}{\mathrm{p}}\right)^{\alpha-\mathrm{m}-1}\text{Ƕ}^3_{\mathrm{qrs}}\left(\frac{\partial^{\mathrm{m}}\mathrm{f(q,r,s,0)}}{\partial \mathrm{t}^{\mathrm{m}}}\right)$$

**Multiplying by $q^{n_1}r^{n_2}s^{n_3}t^{n_4}$**

Let $f(q,r,s,t)$ be a functions such that $\text{Ƕ}^4_{qrst}[f(q,r,s,t)]=F[(h,j,k,l),(m,n,o,p)]$ then $\text{Ƕ}^4_{qrst}[(q^{n_1}r^{n_2}s^{n_3}t^{n_4})f(q,r,s,t)]=$

$$(-1)^{n_1+n_2+n_3+n_4}m^{n_1}n^{n_2}o^{n_3}p^{n_4}\frac{\partial^{n_1+n_2+n_3+n_4}}{\partial h^{n_1}\partial j^{n_2}\partial k^{n_3}\partial l^{n_4}}F[(h,j,k,l),(m,n,o,p)]$$

$\xRightarrow{Proof}$

We know that,

$$\text{Ƕ}^4_{qrst}[f(q,r,s,t)]=\int_0^{\infty}\int_0^{\infty}\int_0^{\infty}\int_0^{\infty}e^{-\left(\frac{hq}{m}+\frac{jr}{n}+\frac{ks}{o}+\frac{lt}{p}\right)}f(q,r,s,t)dqdrdsdt$$

Therefore

$$(-1)^{n_1+n_2+n_3+n_4}\frac{\partial^{n_1+n_2+n_3+n_4}}{\partial h^{n_1}\partial j^{n_2}\partial k^{n_3}\partial l^{n_4}}\text{Ƕ}^4_{qrst}[f(q,r,s,t)]=$$

$$(-1)^{n_1+n_2+n_3+n_4}\frac{\partial^{n_1+n_2+n_3+n_4}}{\partial h^{n_1}\partial j^{n_2}\partial k^{n_3}\partial l^{n_4}}\int_0^{\infty}\int_0^{\infty}\int_0^{\infty}\int_0^{\infty}e^{-\left(\frac{hq}{m}+\frac{jr}{n}+\frac{ks}{o}+\frac{lt}{p}\right)}f(q,r,s,t)dqdrdsdt$$

$$=(-1)^{n_1+n_2+n_3+n_4}\frac{\partial^{n_1+n_2+n_3}}{\partial h^{n_1}\partial j^{n_2}\partial k^{n_3}}\int_0^{\infty}\int_0^{\infty}\int_0^{\infty}e^{-\left(\frac{hq}{m}+\frac{jr}{n}+\frac{ks}{o}\right)}\left[\frac{(-1)^{n_4}}{l^{n_4}}H_t^4(t^{n_4}f(q,r,s,t))\right]dsdrdq$$

$$=\frac{(-1)^{n_1+n_2+n_3}}{l^{n_4}}\frac{\partial^{n_1+n_2+n_3}}{\partial h^{n_1}\partial j^{n_2}\partial k^{n_3}}\int_0^{\infty}\int_0^{\infty}\int_0^{\infty}e^{-\left(\frac{hq}{m}+\frac{jr}{n}+\frac{ks}{o}\right)}[H_t^4(t^{n_4}f(q,r,s,t))]dsdrdq \qquad (20)$$

In the same way integrate above equation with respect to $q,r,and\ s$. We get

For $n_1=n_2=n_3=n_4=1$, we have

$$\text{Ƕ}^4_{qrst}[qrst\,f(q,r,s,t)]=-(mnop)\frac{\partial^4}{\partial h\partial j\partial k\partial k}\text{Ƕ}^4_{qrst}[f(q,r,s,t)]$$

Recall that the Heaviside unit step function $U(q-a,r-b,s-c,t-d)$ (Thakur et al., 2018) is defined by

$$U(q-a,r-b,s-c,t-d)=\begin{cases}1 & q>a,r>b,s>c,t>d\\ 0 & Otherwise\end{cases}$$

**Theorem 3**

Let $f(q,r,s,t)$ be a function such that

$$\text{Ƕ}^4_{qrst}[f(q,r,s,t)]=F[(h,j,k,l),(m,n,o,p)]$$

Then for a constant $a,b,c,d$ we have,

$$\text{Ƕ}^4_{qrst}[f(q-a,r-b,s-c,t-d)U(q-a,r-b,s-c,t-d)]$$
$$=e^{-\left(\frac{hq}{m}+\frac{jr}{n}+\frac{ks}{o}+\frac{lt}{p}\right)}F[(h,j,k,l),(m,n,o,p)]$$

Where $U(q,r,s,t)$ is the Heaviside unit step.

**Proof**

$\Longrightarrow$

$$\text{Ƕ}^4_{qrst}[f(q-a,r-b,s-c,t-d)U(q-a,r-b,s-c,t-d)]=\int_0^{\infty}\int_0^{\infty}\int_0^{\infty}\int_0^{\infty}$$
$$e^{-\left(\frac{hq}{m}+\frac{jr}{n}+\frac{ks}{o}+\frac{lt}{p}\right)}[f(q-a,r-b,s-c,t-d)U(q-a,r-b,s-c,t-d)]dqdrdsdt$$

$$=\int_0^{\infty}\int_0^{\infty}\int_0^{\infty}\int_0^{\infty}e^{-\left(\frac{hq}{m}+\frac{jr}{n}+\frac{ks}{o}+\frac{lt}{p}\right)}[f(q-a,r-b,s-c,t-d)]dqdrdsdt$$

Put $q-a=w_1,r-b=w_2,s-c=w_3,t-d=w_4$

$$= \int_0^\infty \int_0^\infty \int_0^\infty \int_0^\infty e^{-\left[\frac{h(w_1+a)}{m}+\frac{j(w_2+b)}{n}+\frac{k(w_3+c)}{o}+\frac{l(w_4+d)}{p}\right]} f(w_1,w_2,w_3,w_4)\,dw_1dw_2dw_3dw_4$$

$$= e^{-\left(\frac{hq}{m}+\frac{jr}{n}+\frac{ks}{o}+\frac{lt}{p}\right)} F[(h,j,k,l),(m,n,o,p)]$$

## Convolution theorem for the quadruple Shehu Transform

The convolution of the functions $f(q,r,s,t), g(q,r,s,t)$ is denoted by $(f **** g)(q,r,s,t)$ and defined by

$$(f **** g)(q,r,s,t) = \int_0^q \int_0^r \int_0^s \int_0^t f(q-t_1, r-t_2, s-t_3, t-t_4) g(t_1,t_2,t_3,t_4)\,dt_1dt_2dt_3dt_4$$

$$= \int_0^q \int_0^r \int_0^s \int_0^t g(q-t_1, r-t_2, s-t_3, t-t_4) f(t_1,t_2,t_3,t_4)\,dt_1dt_2dt_3dt_4$$

**Theorem 4**

Let $f(q,r,s,t), g(q,r,s,t)$ be of exponential order, such that

$$F[(h,j,k,l),(m,n,o,p)] = \int_0^\infty \int_0^\infty \int_0^\infty \int_0^\infty e^{-\left(\frac{hq}{m}+\frac{jr}{n}+\frac{ks}{o}+\frac{lt}{p}\right)} f(q,r,s,t)\,dqdrdsdt,$$

Is converge, and in addition if

$$G[(h,j,k,l),(m,n,o,p)] = \int_0^\infty \int_0^\infty \int_0^\infty \int_0^\infty e^{-\left(\frac{hq}{m}+\frac{jr}{n}+\frac{ks}{o}+\frac{lt}{p}\right)} g(q,r,s,t)\,dqdrdsdt,$$

Absolutely converge then

$$Ƕ^4_{qrst}[(f **** g)(q,r,s,t)] = Ƕ^4_{qrst}[f(q,r,s,t)]Ƕ^4_{qrst}[g(q,r,s,t)]$$

**Proof**

$\Longrightarrow$ We have,

$$Ƕ^4_{qrst}[(f **** g)(q,r,s,t)] = \int_0^\infty \int_0^\infty \int_0^\infty \int_0^\infty e^{-\left(\frac{hq}{m}+\frac{jr}{n}+\frac{ks}{o}+\frac{lt}{p}\right)} (f **** g)(q,r,s,t)\,dqdrdsdt,$$

$$= \int_0^\infty \int_0^\infty \int_0^\infty \int_0^\infty e^{-\left(\frac{hq}{m}+\frac{jr}{n}+\frac{ks}{o}+\frac{lt}{p}\right)} \left[\int_0^q \int_0^r \int_0^s \int_0^t g(q-t_1, r-t_2, s-t_3, t-t_4) f(t_1,t_2,t_3,t_4)\,dt_1dt_2dt_3dt_4\right] dqdrdsdt$$

Using the Heaviside unit step function, we obtained

$$= \int_0^\infty \int_0^\infty \int_0^\infty \int_0^\infty g(t_1,t_2,t_3,t_4) \left[\int_0^\infty \int_0^\infty \int_0^\infty \int_0^\infty e^{-\left(\frac{hq}{m}+\frac{jr}{n}+\frac{ks}{o}+\frac{lt}{p}\right)} g(q-t_1, r-t_2, s-t_3, t-t_4) f(t_1,t_2,t_3,t_4)\,dqdrdsdt\right] dt_1dt_2dt_3dt_4$$

By using theorem 3,

$$= \int_0^\infty \int_0^\infty \int_0^\infty \int_0^\infty e^{-\left(\frac{ht_1}{m}+\frac{jt_2}{n}+\frac{kt_3}{o}+\frac{lt_4}{p}\right)} F[(h,j,k,l),(m,n,o,p)] g(t_1,t_2,t_3,t_4)\,dt_1dt_2dt_3dt_4$$

$$= F[(h,j,k,l),(m,n,o,p)] \int_0^\infty \int_0^\infty \int_0^\infty \int_0^\infty e^{-\left(\frac{ht_1}{m}+\frac{jt_2}{n}+\frac{kt_3}{o}+\frac{lt_4}{p}\right)} g(t_1,t_2,t_3,t_4)\,dt_1dt_2dt_3dt_4$$

$$= F[(h,j,k,l),(m,n,o,p)]G[(h,j,k,l),(m,n,o,p)]$$

$$= Ƕ^4_{qrst}[f(q,r,s,t)]Ƕ^4_{qrst}[g(q,r,s,t)]$$

## APPLICATIONS

Now we illustrate the quadruple Shehu transform by solving the following examples.

**Example 1**

Consider the following fractional partial differential equation

$$D_s^{\alpha} \text{ш}(q,r,s,t) = \frac{\partial^2 \text{ш}(q,r,s,t)}{\partial q^2}, 0 < \alpha \le 1. \quad (21)$$

With conditions

$$\begin{cases} \text{ш}_q(0,r,s,t) = \sin(r)\, E_{\alpha}(-s^{\alpha}) \\ \text{ш}(0,r,s,t) = 0 \\ \text{ш}(q,r,s,0) = sinq\ sinr\ sins \end{cases} \quad (22)$$

Now by applying quadruple Shehu Transform to equation (21), we have

$$\left(\frac{k}{o}\right)^{\alpha} \widehat{\text{ш}}(q,r,s,t) - \left(\frac{k}{o}\right)^{\alpha-1} \text{ш}(q,r,s,t) = \left(\frac{h}{m}\right)^{2} \widehat{\text{ш}}(q,r,s,t) - \left(\frac{h}{m}\right) \text{Ꞌ}_{rs}^{2} \text{ш}(0,r,s,t) - \text{Ꞌ}_{rs}^{2} \text{ш}_q(0,r,s,t)$$

By using the equation (22), we get

$$\widehat{\text{ш}}[(h,j,k,l),(m,n,o,p)] = \left(\frac{m^2}{h^2+m^2}\right)\left(\frac{n^2}{j^2+n^2}\right)\left(\frac{o^2}{k^2+o^2}\right)\left(\frac{\left(\frac{k}{o}\right)^{\alpha-1}}{\left(\frac{k}{o}\right)^{\alpha}+1}\right) \quad (23)$$

Now, by using the inverse Quadruple Shehu Transform to equation (23), we have

$$\text{ш}(q,r,s,t) = sinq\ sinr\ sins\ E_{\alpha}(-s^{\alpha})$$

**Example 2**

Consider the fourth order homogeneous partial differential equation (Rehman et al., 2014)

$$\frac{\partial^4 \text{ш}(q,r,s,t)}{\partial q \partial r \partial s \partial t} - \text{ш}(q,r,s,t) = 0 \quad (24)$$

With conditions,

$$\text{ш}(q,r,s,0) = e^{q+r+s}, \quad \text{ш}(q,0,s,t) = e^{q+s+t}$$

$$\text{ш}(q,r,0,t) = e^{q+r+t}, \qquad \text{ш}(0,r,s,t) = e^{r+s+t} \quad (25)$$

Now, applying quadruple Shehu Transform to equation (24), we have

$$\left(\frac{hjkl}{mnop} - 1\right) \widehat{\text{ш}}[(h,j,k,l),(m,n,o,p)] - \mathbb{A}[(q,r,s,t)] = 0$$

Where

$$\mathbb{A}[\mathrm{q,r,s,t}] = \left(\frac{hjk}{mno}\right) \text{Ꮒ}_{qrs}^{3}[f(q,r,s,0)] + \left(\frac{hjl}{mnp}\right) \text{Ꮒ}_{qrt}^{3}[f(q,r,0,t)] + \left(\frac{hkl}{mop}\right) \text{Ꮒ}_{qst}^{3}[f(q,0,s,t)] + \left(\frac{jkl}{nop}\right) \text{Ꮒ}_{rst}^{3}[f(0,r,s,t)]$$
$$- \left(\frac{hj}{mn}\right) \text{Ꮒ}_{qr}^{2}[f(q,r,0,0)] - \left(\frac{hk}{mo}\right) \text{Ꮒ}_{qs}^{2}[f(q,0,s,0)] - \left(\frac{hl}{mp}\right) \text{Ꮒ}_{qt}^{2}[f(q,0,0,t)] - \left(\frac{jk}{no}\right) \text{Ꮒ}_{rs}^{2}[f(0,r,s,0)]$$
$$- \left(\frac{jl}{np}\right) \text{Ꮒ}_{rt}^{2}[f(0,r,0,t)] - \left(\frac{kl}{op}\right) \text{Ꮒ}_{st}^{2}[f(0,0,s,t)] + \left(\frac{h}{m}\right) \text{Ꮒ}_{q}[f(q,0,0,0)] + \left(\frac{j}{n}\right) \text{Ꮒ}_{r}[f(0,r,0,0)]$$
$$+ \left(\frac{k}{o}\right) \text{Ꮒ}_{s}[f(0,0,s,0)] + \left(\frac{l}{p}\right) \text{Ꮒ}_{t}[f(0,0,0,t)] - f(0,0,0,0)$$

$$= \left(\frac{m}{h-m}\right)\left(\frac{n}{j-n}\right)\left(\frac{o}{k-o}\right)\left(\frac{p}{l-p}\right)\left(\frac{hjkl-mnop}{mnop}\right)$$

By using equation (25), we get

$$\widehat{\text{ɯ}}[(h,j,k,l),(m,n,o,p)] = \left(\frac{m}{h-m}\right)\left(\frac{n}{j-n}\right)\left(\frac{o}{k-o}\right)\left(\frac{p}{l-p}\right) \tag{26}$$

By taking inverse quadruple Shehu Transform of equation ( 26)

$$\text{ɯ}(q,r,s,t) = \text{Ꮒ}_{qrst}^{-4}\left[\left(\frac{m}{h-m}\right)\left(\frac{n}{j-n}\right)\left(\frac{o}{k-o}\right)\left(\frac{p}{l-p}\right)\right] = e^{q+r+s+t}$$

**Example 3**

Consider the fourth order non-homogeneous partial differential equation (Rehman et al., 2014)

$$\frac{\partial^4 \text{ɯ}(q,r,s,t)}{\partial q \partial r \partial s \partial t} + \text{ɯ}(q,r,s,t) = 5e^{-2q+r-2s+t} \tag{27}$$

With conditions,

$$\text{ɯ}(q,r,s,0) = e^{-2q+r-2s}, \quad \text{ɯ}(q,0,s,t) = e^{-2q-2s+t}$$

$$\text{ɯ}(q,r,0,t) = e^{-2q+r+t}, \qquad \text{ɯ}(0,r,s,t) = e^{r-2s+t} \tag{28}$$

Now, applying quadruple Shehu Transform to equation ( 28), we have

$$\left(\frac{hjkl}{mnop}+1\right)\widehat{\text{ɯ}}[(h,j,k,l),(m,n,o,p)] - \mathbb{A}[(q,r,s,t)] = 5\left(\frac{m}{h+2m}\right)\left(\frac{n}{j-n}\right)\left(\frac{o}{k+2o}\right)\left(\frac{p}{l-p}\right)$$

$$\left(\frac{hjkl}{mnop}+1\right)\widehat{\text{ɯ}}[(h,j,k,l),(m,n,o,p)] = \mathbb{A}[(q,r,s,t)] + 5\left(\frac{m}{h+2m}\right)\left(\frac{n}{j-n}\right)\left(\frac{o}{k+2o}\right)\left(\frac{p}{l-p}\right)$$

Where

$$\mathbb{A}[\mathrm{q,r,s,t}] = \left(\frac{hjk}{mno}\right) \text{Ꮒ}_{qrs}^{3}[\text{ɯ}(q,r,s,0)] + \left(\frac{hjl}{mnp}\right) \text{Ꮒ}_{qrt}^{3}[\text{ɯ}(q,r,0,t)] + \left(\frac{hkl}{mop}\right) \text{Ꮒ}_{qst}^{3}[\text{ɯ}(q,0,s,t)]$$
$$+ \left(\frac{jkl}{nop}\right) \text{Ꮒ}_{rst}^{3}[\text{ɯ}(0,r,s,t)] - \left(\frac{hj}{mn}\right) \text{Ꮒ}_{qr}^{2}[\text{ɯ}(q,r,0,0)] - \left(\frac{hk}{mo}\right) \text{Ꮒ}_{qs}^{2}[\text{ɯ}(q,0,s,0)]$$
$$- \left(\frac{hl}{mp}\right) \text{Ꮒ}_{qt}^{2}[\text{ɯ}(q,0,0,t)] - \left(\frac{jk}{no}\right) \text{Ꮒ}_{rs}^{2}[\text{ɯ}(0,r,s,0)] - \left(\frac{jl}{np}\right) \text{Ꮒ}_{rt}^{2}[\text{ɯ}(0,r,0,t)]$$
$$- \left(\frac{kl}{op}\right) \text{Ꮒ}_{st}^{2}[\text{ɯ}(0,0,s,t)] + \left(\frac{h}{m}\right) \text{Ꮒ}_{q}[\text{ɯ}(q,0,0,0)] + \left(\frac{j}{n}\right) \text{Ꮒ}_{r}[\text{ɯ}(0,r,0,0)]$$
$$+ \left(\frac{k}{o}\right) \text{Ꮒ}_{s}[\text{ɯ}(0,0,s,0)] + \left(\frac{l}{p}\right) \text{Ꮒ}_{t}[\text{ɯ}(0,0,0,t)] - \text{ɯ}(0,0,0,0)$$

$$= \left(\frac{hjkl - 4mnop}{(h+2m)(j-n)(k+2o)(l-p)}\right)$$

By using equation (28), we get

$$\widehat{\mathrm{w}}[(h,j,k,l),(m,n,o,p)] = \left(\frac{m}{h+2m}\right)\left(\frac{n}{j-n}\right)\left(\frac{o}{k+2o}\right)\left(\frac{p}{l-p}\right) \quad (29)$$

By taking inverse quadruple Shehu Transform of equation (29)

$$\mathrm{w}(q,r,s,t) = \mathrm{H}_{qrst}^{-4}\left[\left(\frac{m}{h+2m}\right)\left(\frac{n}{j-n}\right)\left(\frac{o}{k+2o}\right)\left(\frac{p}{l-p}\right)\right] = e^{-2q+r-2s+t}$$

## Conclusion

In this work, we introduced a new type of integral transform called quadruple Shehu transform. First, we introduced existence and uniqueness of new transform (i.e. quadruple Shehu transform) and their proofs. We also, applied the new transform on some functions. To illustrate the efficiency of quadruple Shehu transform, we applied this transform on some examples. We got that quadruple Shehu transform is a more efficient method for solving homogeneous and non-homogeneous partial differential equation also fractional partial differential equations. Now this new transform will extend to known results on a fourth Laplace transform and Fourth Aboodh transform to our results on a quadruple Shehu transform.